\RequirePackage{ifpdf}
\ifpdf 
\documentclass[pdftex]{sigma}
\else
\documentclass{sigma}
\fi

\DeclareMathOperator{\domain}{Dom}

\newcommand{\N}{\mathbb{N}^d}
\newcommand{\R}{\mathbb{R}^d}
\newcommand{\C}{\mathbb{C}}

\begin{document}

\allowdisplaybreaks

\renewcommand{\thefootnote}{$\star$}

\renewcommand{\PaperNumber}{016}

\FirstPageHeading

\ShortArticleName{Imaginary Powers of the Dunkl Harmonic Oscillator}

\ArticleName{Imaginary Powers of the Dunkl Harmonic Oscillator\footnote{This paper is a contribution to the Special
Issue on Dunkl Operators and Related Topics. The full collection
is available at
\href{http://www.emis.de/journals/SIGMA/Dunkl_operators.html}{http://www.emis.de/journals/SIGMA/Dunkl\_{}operators.html}}}

\Author{Adam NOWAK and Krzysztof STEMPAK}

\AuthorNameForHeading{A. Nowak and K. Stempak}

\Address{Instytut Matematyki i Informatyki, Politechnika Wroc\l{}awska,   \\
      Wyb. Wyspia\'nskiego 27, 50--370 Wroc\l{}aw, Poland}

\Email{\href{mailto:Adam.Nowak@pwr.wroc.pl}{Adam.Nowak@pwr.wroc.pl}, \href{mailto:Krzysztof.Stempak@pwr.wroc.pl}{Krzysztof.Stempak@pwr.wroc.pl}}

\URLaddress{\url{http://www.im.pwr.wroc.pl/~anowak/}, \url{http://www.im.pwr.wroc.pl/~stempak/}}

\ArticleDates{Received October 14, 2008, in f\/inal form February 08,
2009; Published online February 11, 2009}

\Abstract{In this paper we continue the study of spectral properties of the Dunkl
harmo\-nic~oscillator in the context of a f\/inite ref\/lection group on $\R$
isomorphic to $\mathbb{Z}^d_2$. We prove that imaginary powers of this operator are
bounded on $L^p$, $1<p<\infty$, and from $L^1$ into weak~$L^1$.}


\Classification{42C10; 42C20}

\section{Introduction} \label{sec:intro}

In \cite{NS2} the authors def\/ined and investigated a system of Riesz transforms related to the Dunkl
harmonic oscillator $\mathcal{L}_{k}$. The present article continues the study of spectral properties
of operators associated with $\mathcal{L}_{k}$ by considering the imaginary powers
$\mathcal{L}_{k}^{-i\gamma}$, $\gamma\in\mathbb{R}$.
Our objective is to study $L^p$ mapping properties of the operators $\mathcal{L}_{k}^{-i\gamma}$,
and the principal tool is the general Calder\'on--Zygmund operator theory.
The main result we get (Theorem \ref{thm:main})
partially extends the result obtained recently by Stempak and Torrea \cite[Theorem 4.3]{ST2}
and corresponding to the trivial multiplicity function $k \equiv 0$.
Imaginary powers of the Euclidean Laplacian were investigated much earlier by Muckenhoupt \cite{M}.

Let us brief\/ly describe the framework of the Dunkl theory of dif\/ferential-dif\/ference operators on~$\R$ related to f\/inite ref\/lection groups. Given such a group $G\subset O(\R)$ and a $G$-invariant
nonnegative multiplicity function $k\colon R\to[0,\infty)$ on a root system $R\subset \R$ associated
with the ref\/lections of~$G$, the Dunkl dif\/ferential-dif\/ference operators $T_j^{k}$, $j=1,\ldots,d$,
are def\/ined by
\begin{gather*}
T_j^{k}f(x)=\partial_j f(x) +\sum_{\beta\in R_+}k(\beta)\beta_j\frac{f(x)-
	f(\sigma_\beta x)}{\langle\beta,x\rangle}, \qquad f\in C^1(\R);
\end{gather*}
here $\partial_j$ is the $j$th partial derivative, $\langle\cdot,\cdot\rangle$ denotes the standard
inner product in~$\R$, $R_+$ is a~f\/ixed positive subsystem of $R$, and $\sigma_\beta$ denotes the
ref\/lection in the hyperplane orthogonal to~$\beta$. The Dunkl operators $T_j^{k}$, $j=1,\ldots,d$,
form a commuting system (this is an important feature, see \cite{Du1})
of the f\/irst order dif\/ferential-dif\/ference
operators, and reduce to $\partial_j$, $j=1,\ldots,d$, when $k\equiv0$. Moreover,
$T_j^{k}$ are homogeneous of degree $-1$ on~$\mathcal{P}$, the space of all polynomials in~$\R$. This means that $T_j^{k}\mathcal{P}_{m}\subset\mathcal{P}_{m-1}$, where
$m\in\mathbb{N}=\{0,1,\ldots\}$ and $\mathcal{P}_{m}$ denotes the subspace of~$\mathcal{P}$
consisting of polynomials of total degree $m$ (by convention, $\mathcal{P}_{-1}$ consists only
of the null function).

In Dunkl's theory the operator, see \cite{Du2},
\begin{gather*}
\Delta_{k}=\sum_{j=1}^d (T_j^{k})^2
\end{gather*}
plays the role of the Euclidean Laplacian (notice that $\Delta$ comes into play when $k\equiv0$).
It is homogeneous of degree $-2$ on $\mathcal{P}$ and symmetric in $L^2(\R, w_{k})$, where
\begin{gather*}
w_{k}(x)=\prod_{\beta\in R_+}|\langle\beta,x\rangle|^{2k(\beta)},
\end{gather*}
if considered initially on $C^\infty_{c}(\R)$. Note that $w_{k}$ is $G$-invariant.

The study of the operator
\begin{gather*}
L_{k}=-\Delta_{k}+\|x\|^2
\end{gather*}
was initiated by R\"osler \cite{R1,R2}. It occurs that $L_{k}$ (or rather its self-adjoint extension
$\mathcal{L}_{k}$) has a discrete spectrum and the corresponding eigenfunctions are the
generalized Hermite functions def\/ined and investigated by R\"osler \cite{R1}.
Due to the form of $L_k$, it is reasonable to call it the {\sl Dunkl harmonic oscillator}.
In fact $L_k$ becomes the classic harmonic oscillator $-\Delta+\|x\|^2$ when $k\equiv0$.

The results of the present paper are naturally related to the authors' articles
\cite {NS1,NS2}. In what follows we will use the notation introduced there and
invoke certain arguments from \cite {NS1}. For basic facts concerning Dunkl's theory
we refer the reader to the excellent survey article by R\"osler~\cite{R3}.

Throughout the paper we use a fairly standard notation.
Given a multi-index $n\in\N$, we write $|n|=n_1+\dots+n_d$
and, for $x,y\in\R$,  $xy = (x_1 y_1,\ldots,x_d y_d)$,
$x^n=x_1^{n_1}\cdot\cdots\cdot x_d^{n_d}$ (and similarly $x^{\alpha}$ for $x\in\R_+$ and $\alpha\in\R$);
$\|x\|$ denotes the Euclidean norm of $x\in\R$, and $e_j$ is the $j$th coordinate vector in $\R$.
Given $x\in \R$ and $r>0$, $B(x,r)$ is the Euclidean ball in $\R$ centered at $x$ and of radius $r$.
For a nonnegative weight function $w$ on $\R$, by $L^p(\R,w)$, $1\le p<\infty$, we denote the
usual Lebesgue spaces related to the measure $dw(x)=w(x)dx$
(in the sequel we will often abuse slightly the notation and use the same symbol $w$ to
denote the measure induced by a density $w$).
Writing $X\lesssim Y$ indicates that $X\leq CY$ with a positive constant $C$
independent of signif\/icant quantities. We shall write $X \simeq Y$ when $X \lesssim Y$ and $Y \lesssim X.$

\section{Preliminaries} \label{sec:prel}

In the setting of general Dunkl's theory R\"osler \cite{R1} constructed systems of
naturally associated multivariable generalized Hermite polynomials and Hermite functions.
The system of generalized Hermite polynomials $\{H_n^{k} : n\in\N\}$
is orthogonal and complete in $L^2(\R, e^{-\|\cdot\|^2}w_{k})$, while the system $\{h_n^{k} : n\in\N\}$
of generalized Hermite functions
\begin{gather*}
h_n^{k}(x)=\big({2^{|n|}c_{k}}\big)^{-1/2}\exp(-\|x\|^2/2)H_n^{k}(x), \qquad x\in \R,\quad n\in \N,
\end{gather*}
is an orthonormal basis in $L^2(\R, w_{k})$, cf.~\cite[Corollary 3.5 (ii)]{R1};
here the normalizing constant~$c_k$ equals to $\int_{\R}\exp(-\|x\|^2)w_{k}(x)\,dx$.
Moreover,
$h_n^k$ are eigenfunctions of  $L_{k}$,
\begin{gather*}
L_{k}h_n^{k}=(2|n|+2\tau+d)h_n^{k},
\end{gather*}
where $\tau=\sum_{\beta\in R_+}k(\beta)$. For $k\equiv0$,
$h_n^{0}$ are the usual multi-dimensional Hermite functions, see for instance \cite{ST1} or \cite{ST2}.

Let $\langle\cdot,\cdot\rangle_{k}$ be the canonical inner product in $L^2(\mathbb{R}^d, w_{k})$.
The operator
\begin{gather*}
  \mathcal{L}_{k} f=\sum_{n\in \N} (2|n|+2\tau+d)\langle f,h_n^{k}
  \rangle_{{k}} \, h_n^{k},
\end{gather*}
def\/ined on the domain
\begin{gather*}
  \domain(\mathcal{L}_{k})=\Big\{f\in L^2(\R, w_{{k}}):
  \sum_{n\in \N}\big|(2|n|+2\tau+d)\langle f,h_n^k
  \rangle_{k}\big|^2<\infty\Big\},
\end{gather*}
is a self-adjoint extension of $L_{k}$ considered on $C^{\infty}_{c}(\R)$
as the natural domain (the inclusion $C^\infty_{c}(\R)\subset
\domain(\mathcal{L}_{k})$ may be easily verif\/ied).
The spectrum of $\mathcal{L}_k$ is the discrete set $\{2m+2\tau+d:m\in \mathbb{N}\}$, and
the spectral decomposition of $\mathcal{L}_k$ is
\begin{gather*}
\mathcal{L}_{k} f=\sum_{m=0}^\infty(2m+2\tau+d)\mathcal{P}_m^k f, \qquad
f\in\domain(\mathcal{L}_{k}),
\end{gather*}
where the spectral projections are
\begin{gather*}
\mathcal{P}_m^{k} f=\sum_{|n|=m}\langle f,h_n^{k}
  \rangle_{{k}} \,h_n^{k}.
\end{gather*}
By Parseval's identity, for each $\gamma\in\mathbb{R}$ the operator
\begin{gather*}
\mathcal{L}_{k}^{-i\gamma}f=\sum_{m=0}^\infty (2m+2\tau+d)^{-i\gamma}\, \mathcal{P}_m^{k} f
\end{gather*}
is an isometry on $L^2(\R, w_{k})$.

Consider the f\/inite ref\/lection group generated by $\sigma_j$, $j=1,\ldots,d$,
\begin{gather*}
\sigma_{j}(x_1,\ldots,x_j,\ldots,x_d)=(x_1,\ldots,-x_j,\ldots,x_d),
\end{gather*}
and isomorphic to $\mathbb{Z}_2^d = \{0,1\}^d$.
The ref\/lection $\sigma_j$ is in the hyperplane orthogonal to $e_j$.
Thus $R=\{\pm \sqrt{2}e_j: j=1,\ldots,d\}$, $R_+=\{\sqrt{2}e_j: j=1,\ldots,d\}$,
and for a nonnegative multiplicity function $k\colon R\to [0,\infty)$ which is
$\mathbb{Z}^d_2$-invariant only values of $k$ on $R_+$ are essential. Hence we may think
$k=(\alpha_1+1\slash 2,\ldots,\alpha_d+1\slash 2)$,
$\alpha_j\ge-1/2$. We write $\alpha_j+1\slash 2$ in place of seemingly more appropriate $\alpha_j$
since, for the sake of clarity, it is convenient for us to stick to the notation used
in \cite {NS1} and \cite{NS2}.

In what follows the symbols $T_j^\alpha$,  $\Delta_\alpha$, $w_\alpha$, $L_\alpha$, $\mathcal L_\alpha$,
$h_n^\alpha$, and so on, denote the objects introduced earlier
and related to the present $\mathbb{Z}_2^d$ group setting.
Thus the Dunkl dif\/ferential-dif\/ference operators are now given by
\begin{gather*}
T_j^{\alpha}f(x)=\partial_j f(x) +(\alpha_j+1\slash 2)\frac{f(x)-f(\sigma_j x)}{x_j},
\qquad f\in C^1(\R),
\end{gather*}
and the explicit formula for the Dunkl Laplacian is
\begin{gather*}
\Delta_\alpha f(x)=
\sum_{j=1}^d\left(\frac{\partial^2f}{\partial x^2_j}(x)+\frac{2\alpha_j+1}{x_j}\frac{\partial
f}{\partial x_j}(x)-(\alpha_j+1\slash 2)\frac{f(x)-f(\sigma_jx)}{x_j^2}\right).
\end{gather*}
The corresponding weight $w_\alpha$ has the form
\begin{gather*}
w_\alpha(x)= \prod_{j=1}^d|x_j|^{2\alpha_j+1} \simeq
\prod_{\beta \in R_+} |\langle \beta, x \rangle_{\alpha}|^{2k(\beta)}, \qquad x\in \R.
\end{gather*}
Given $\alpha \in [-1/2,\infty)^{d}$,
the associated generalized Hermite functions are tensor products
\begin{gather*}
h_{n}^{\alpha}(x) = h _{n_1}^{\alpha _1}(x_1) \cdot \cdots \cdot
h_{n_d}^{\alpha _d}(x_d), \qquad x = (x_1, \ldots ,x_d)\in \R,\quad n=(n_1,\ldots,n_d) \in \N,
\end{gather*}
where $h_{n_i}^{\alpha _i}$ are the one-dimensional functions (see Rosenblum \cite{Ros})
\begin{gather*}
h_{2n_i}^{\alpha_i}(x_i) =d_{2n_i,\alpha_i}e^{-x_i^2/2}L_{n_i}^{\alpha_i}\big(x_i^2\big),\\
h_{2n_i+1}^{\alpha _i}(x_i) =d_{2n_i+1,\alpha_i}e^{-x_i^2/2}x_iL_{n_i}^{\alpha_i+1}\big(x_i^2\big);
\end{gather*}
here $L_{n_i}^{\alpha_i}$ denotes the Laguerre polynomial of degree $n_i$ and order $\alpha_i$,
cf.~\cite[p.~76]{Leb}, and
\begin{gather*}
d_{2n_i,\alpha_i}=(-1)^{n_i}\bigg(\frac{\Gamma(n_i+1)}{\Gamma(n_i+\alpha_i +1)}\bigg)^{1/2}, \qquad
d_{2n_i+1,\alpha_i}=(-1)^{n_i}\bigg(\frac{\Gamma(n_i+1)}{\Gamma(n_i+\alpha_i +2)}\bigg)^{1/2}.
\end{gather*}
For $\alpha=(-1/2,\ldots,-1/2)$ we obtain the usual Hermite functions.
The system $\{h_{n}^{\alpha} : n \in \N \}$ is an orthonormal basis in $L^2(\R,w_{\alpha})$ and
\begin{gather*}
L_\alpha h_n^{\alpha}=(2|n|+2|\alpha |+2d)h_n^{\alpha},
\end{gather*}
where by $|\alpha|$ we denote $|\alpha|=\alpha_1+\dots+\alpha_d$ (thus $|\alpha|$ may be negative).

The semigroup $T_t^\alpha = \exp(-t\mathcal{L}_{\alpha})$, $t \ge 0$, generated by
$\mathcal{L}_{\alpha}$ is a strongly continuous semigroup of contractions on
$L^2(\R, w_{\alpha})$. By the spectral theorem,
\begin{gather*} \label{sr_t}
T_t^\alpha f=\sum_{m=0}^\infty e^{-t(2m+2|\alpha|+2d)}\mathcal{P}^\alpha_mf,
\qquad f\in L^2(\R, w_{\alpha}).
\end{gather*}
The integral representation of $T_t^\alpha$ on $L^2(\R,w_{\alpha})$ is
\begin{gather*}
  T_t^\alpha f(x)=\int_{\R} G_t^\alpha(x,y)f(y)\,dw_{\alpha}(y),
  \qquad x\in \R, \quad t>0,
\end{gather*}
where the heat kernel $\{G_t^\alpha\}_{t>0}$ is given by
\begin{gather} \label{jc}
G^\alpha_t(x,y)=\sum_{m=0}^\infty e^{-t(2m+2|\alpha|+2d)} \sum_{|n|=m}
h_n^\alpha(x)h_n^\alpha(y).
\end{gather}

In dimension one, for $\alpha\ge-1/2$ it is known (see for instance
\cite[Theorem 3.12]{R1} and \cite[p.~523]{R1}) that
\begin{gather*}
G^\alpha_t(x,y)=\frac{1}{2\sinh 2t}\exp\left({-\frac{1}{2} \coth(2t)\big(x^{2}+y^{2}\big)}\right)
\left[\frac{I_{\alpha}\left(\frac{x y}{\sinh 2t}\right)}
{(x y)^{\alpha}}+xy \frac{I_{\alpha+1}\left(\frac{x y}{\sinh 2t}\right)}{(x y)^{\alpha+1}}\right],
\end{gather*}
with $I_\nu$ being the modif\/ied Bessel function of the f\/irst kind and order $\nu$,
\begin{gather*}
I_{\nu}(z) = \sum_{k=0}^{\infty} \frac{(z\slash 2)^{\nu+2k}}{\Gamma(k+1)\Gamma(k+\nu+1)}.
\end{gather*}
Here we consider the function $z \mapsto z^{\nu}$, and thus also the Bessel function $I_{\nu}(z)$,
as an analytic function def\/ined on $\C \backslash \{ix : x \le 0\}$ (usually $I_{\nu}$ is considered
as a function on $\C$ cut along the half-line $(-\infty,0]$).
Note that $I_{\nu}$, as a function on $\mathbb{R}_{+}$,
is real, positive and smooth for any $\nu > -1$, see \cite[Chapter 5]{Leb}.

Therefore, in $d$ dimensions,
\begin{gather*}
G^\alpha_t(x,y)=\sum_{\varepsilon\in \mathbb{Z}_2^d}G^{\alpha,\varepsilon}_t(x,y),
\end{gather*}
where the component kernels are
\begin{gather*}
G^{\alpha,\varepsilon}_t(x,y)=\frac{1}{(2\sinh 2t)^{d}}\exp\left({-\frac{1}{2}
 \coth(2t)\big(\|x\|^{2}+\|y\|^{2}\big)}\right) \prod^{d}_{i=1} (x_i y_i)^{\varepsilon_i}
  \frac{I_{\alpha_i+\varepsilon_i}\left(\frac{x_i y_i}{\sinh 2t}\right)}{(x_i
   y_i)^{\alpha_i+\varepsilon_i}}.
\end{gather*}
Note that $G^{\alpha,\varepsilon}_t(x,y)$ is given by the series \eqref{jc}, with the summation
in $n$ restricted to the set of multi-indices
\begin{gather*}
\mathcal{N}_{\varepsilon} = \big\{ n \in \mathbb{N}^d : n_i \; \textrm{is even if} \; \varepsilon_i=0
	\;\textrm{or}\; n_i \; \textrm{is odd if}\; \varepsilon_i=1,\; i=1,\ldots,d\big\}.
\end{gather*}
To verify this fact it is enough to restrict to the one-dimensional case and then
use the Hille--Hardy formula, cf.~\cite[(4.17.6)]{Leb}.

In the sequel we will make use of the following technical result concerning $G^{\alpha,\varepsilon}_t(x,y)$.
The corresponding proof is given at the end of Section \ref{sec:kernel}.
\begin{lemma} \label{der_est}
Let $\alpha \in [-1\slash 2,\infty)^d$ and let $\varepsilon \in \mathbb{Z}_2^d$. Then, with
$x,y\in \mathbb{R}^d_+$ fixed, $x \neq y$, the kernel $G_t^{\alpha,\varepsilon}(x,y)$ decays rapidly
when either $t \to 0^+$ or $t \to \infty$.
Further, given any disjoint compact sets $E,F \subset \mathbb{R}^d_+$, we have
\begin{gather} \label{g_e}
\int_0^{\infty} \big| \partial_t G_t^{\alpha,\varepsilon}(x,y)\big| \, dt \lesssim 1,
\end{gather}
uniformly in $x \in E$ and $y \in F$.
\end{lemma}

We end this section with pointing
out that there is a general background for the facts conside\-red here for an arbitrary ref\/lection
group, see \cite{R3} for a comprehensive account.
In particular, the heat (or Mehler) kernel \eqref{jc} has always a closed form involving the
so-called Dunkl kernel, and is always strictly positive.
This implies that the corresponding semigroup is contractive on $L^{\infty}(\mathbb{R}^d,w_k)$,
and as its generator is self-adjoint and positive in $L^2(\mathbb{R}^d,w_k)$, the semigroup
is also contractive on the latter space. Hence, by duality and interpolation, it is in fact
contractive on all $L^p(\mathbb{R}^d,w_k)$, $1\le p \le \infty$.

\section{Main result} \label{sec:main}
From now on we assume $\gamma\in\mathbb{R}$, $\gamma\neq0$, to be f\/ixed.
Recall that the operator $\mathcal{L}_{\alpha}^{-i\gamma}$ is given on $L^2(\mathbb{R}^d,w_{\alpha})$
by the spectral series,
\begin{gather*}
\mathcal{L}_{\alpha}^{-i\gamma} f = \sum_{n \in \mathbb{N}^d} (2|n|+2|\alpha|+2d)^{-i\gamma}
	\langle f,h_n^{\alpha} \rangle_{\alpha} \,h_n^{\alpha}.
\end{gather*}
Our main result concerns mapping properties of $\mathcal{L}_{\alpha}^{-i\gamma} f$ on $L^p$ spaces.
\begin{theorem}  \label{thm:main}
Assume $\alpha \in [-1\slash 2,\infty)^d$.
Then $\mathcal{L}^{-i\gamma}_{\alpha}$, defined initially on $L^2(\mathbb{R}^d,w_{\alpha})$,
extends uniquely to a bounded operator on $L^p(\R,w_{\alpha})$,
$1<p<\infty$, and to a bounded operator from $L^1(\R,w_{\alpha})$ to $L^{1,\infty}(\R,w_{\alpha})$.
\end{theorem}
The proof we give relies on splitting $\mathcal{L}_{\alpha}^{-i\gamma}$ in $L^2(\R,w_{\alpha})$
into a f\/inite number of suitable $L^2$-bounded operators and then treating each of the operators
separately. More precisely, we decompose
\begin{gather*}
\mathcal{L}_{\alpha}^{-i\gamma}
= \sum_{\varepsilon \in \mathbb{Z}_2^d} \mathcal{L}_{\alpha,\varepsilon}^{-i\gamma},
\end{gather*}
where (with the set $\mathcal{N}_{\varepsilon}$ introduced in the previous section)
\begin{gather*}
\mathcal{L}_{\alpha,\varepsilon}^{-i\gamma} f =
	\sum_{n \in \mathcal{N}_{\varepsilon}} (2|n|+2|\alpha|+2d)^{-i\gamma}
	\langle f,h_n^{\alpha} \rangle_{\alpha} \,h_n^{\alpha}, \qquad f \in L^2(\R,w_{\alpha}).
\end{gather*}
Clearly, each $\mathcal{L}_{\alpha,\varepsilon}^{-i\gamma}$ is a contraction in $L^2(\R,w_{\alpha})$.

It is now convenient to introduce the following terminology:
given $\varepsilon \in \mathbb{Z}_2^d$, we say that a~function $f$ on $\R$ is $\varepsilon$-\emph{symmetric}
if for each $i=1,\ldots,d$, $f$ is either even or odd with respect to the $i$th coordinate
according to whether $\varepsilon_i=0$ or $\varepsilon_i=1$, respectively.
Thus $f$ is $\varepsilon$-symmetric if and only if $f\circ \sigma_i = (-1)^{\varepsilon_i} f$,
$i=1,\ldots,d$. Any function $f$ on $\R$ can be split uniquely into a~sum of $\varepsilon$-symmetric
functions $f_{\varepsilon}$,
\begin{gather*}
f = \sum_{\varepsilon \in \mathbb{Z}_2^d} f_{\varepsilon}, \qquad
	f_{\varepsilon}(x) = \frac{1}{2^d} \sum_{\eta \in \{-1,1\}^d} \eta^{\varepsilon} f(\eta x).
\end{gather*}
For $f \in L^2(\R,w_{\alpha})$ this splitting is orthogonal in $L^2(\R,w_{\alpha})$.
Finally, notice that $h_n^{\alpha}$ is $\varepsilon$-symmetric if and only if
$n \in \mathcal{N}_{\varepsilon}$.
Consequently, $\mathcal{L}_{\alpha,\varepsilon}^{-i\gamma}$ is invariant on the subspace of
$L^2(\R,w_{\alpha})$ of $\varepsilon$-symmetric functions and vanishes on the orthogonal
complement of that subspace.

Observe that in order to prove Theorem \ref{thm:main} it is suf\/f\/icient to show the analogous result
for each $\mathcal{L}_{\alpha,\varepsilon}^{-i\gamma}$. Moreover, since
\begin{gather*}
\mathcal{L}_{\alpha}^{-i\gamma} f = \sum_{\varepsilon \in \mathbb{Z}_2^d}
	\mathcal{L}_{\alpha,\varepsilon}^{-i\gamma} f =
	\sum_{\varepsilon \in \mathbb{Z}_2^d} \mathcal{L}_{\alpha,\varepsilon}^{-i\gamma} f_{\varepsilon}
\end{gather*}
and since for a f\/ixed $1\le p < \infty$ (recall that $w_{\alpha}(\xi x)=w_{\alpha}(x)$,
$\xi \in \{-1,1\}^d$)
\begin{gather*}
\|f\|_{L^p(\R,w_{\alpha})} \simeq \sum_{\varepsilon \in \mathbb{Z}_2^d}
	\|f_{\varepsilon}\|_{L^p(\mathbb{R}^d_+,w_{\alpha}^+)},
\end{gather*}
it is enough to restrict the situation to the space $(\mathbb{R}_+^d,w_{\alpha}^+)$, where
$w_{\alpha}^+$ is the restriction of $w_{\alpha}$ to~$\mathbb{R}^d_+$.
Thus we are reduced to considering the operators
\begin{gather} \label{brs}
\mathcal{L}_{\alpha,\varepsilon,+}^{-i\gamma} f =
	\sum_{n \in \mathcal{N}_{\varepsilon}} (2|n|+2|\alpha|+2d)^{-i\gamma}
	\langle f,h_n^{\alpha} \rangle_{L^2(\mathbb{R}^d_+,w_{\alpha}^+)} \,h_n^{\alpha},
		\qquad f \in L^2(\mathbb{R}^d_+,w^+_{\alpha}),
\end{gather}
which are bounded on $L^2(\mathbb{R}^d_+,w_{\alpha}^+)$
since the system $\{2^{d\slash 2}h_n^{\alpha} : n \in \mathcal{N}_{\varepsilon}\}$
is orthonormal in $L^2(\mathbb{R}^d_+,w_{\alpha}^+)$.
Now, Theorem \ref{thm:main} will be justif\/ied once we prove the following.
\begin{lemma} \label{lem:main}
Assume that $\alpha \in [-1\slash 2,\infty)^d$ and $\varepsilon \in \mathbb{Z}_2^d$. Then
$\mathcal{L}^{-i\gamma}_{\alpha,\varepsilon,+}$, defined initially on $L^2(\mathbb{R}_+^d,w^+_{\alpha})$,
extends uniquely to a bounded operator on $L^p(\mathbb{R}^d_+,w^+_{\alpha})$,
$1<p<\infty$, and to a~bounded operator
from $L^1(\mathbb{R}^d_+,w^+_{\alpha})$ to $L^{1,\infty}(\mathbb{R}^d_+,w^+_{\alpha})$.
\end{lemma}
The proof of Lemma \ref{lem:main} will be furnished by means of the general Calder\'on--Zygmund
theory. In fact, we shall show that each $\mathcal{L}^{-i\gamma}_{\alpha,\varepsilon,+}$
is a Calder\'on--Zygmund operator in the sense of the space of homogeneous type
$(\mathbb{R}_+^d,w_{\alpha}^+,\|\cdot\|)$.
It is well known that the classical Calder\'on--Zygmund theory works,
with appropriate adjustments, when the underlying space is of homogeneous type.
Thus we shall use properly adjusted facts from the classic
Calder\'on--Zygmund theory (presented, for instance, in \cite{Duo})
in the setting of the space $(\R_+,w^+_{\alpha},\|\cdot\|)$ without further comments.

A formal computation based on the formula
\begin{gather*}
\lambda^{-i \gamma} = \frac{1}{\Gamma(i\gamma)} \int_0^{\infty} e^{-t \lambda} t^{i\gamma -1} \, dt,
	\qquad \lambda >0,
\end{gather*}
suggests that $\mathcal L^{-i\gamma}_{\alpha,\varepsilon,+}$ should be associated with the kernel
\begin{gather}\label{kernel}
K_\gamma^{\alpha,\varepsilon}(x,y)=\frac1{\Gamma(i\gamma)}\int_0^\infty
	G_t^{\alpha,\varepsilon}(x,y)t^{i\gamma-1}\,dt, \qquad x,y \in \mathbb{R}^d_+
\end{gather}
(note that for $x \neq y$ the last integral is absolutely convergent due to the decay of
$G_t^{\alpha,\varepsilon}(x,y)$ at $t\to0^+$ and $t\to\infty$, see Lemma \ref{der_est}).
The next result shows that this is indeed the case, at least in the Calder\'on--Zygmund theory sense.
\begin{proposition}  \label{sto}
Let $\alpha \in [-1\slash 2, \infty)^d$ and $\varepsilon \in \mathbb{Z}_2^d$. Then for
$f,g\in C^\infty_{c}({\R_+})$ with disjoint supports
\begin{gather}
\langle \mathcal L^{-i\gamma}_{\alpha,\varepsilon,+}f,g\rangle_{L^2(\R_+,w^+_{\alpha})}=
\int_{{\R_+}}\int_{{\R_+}}
K_\gamma^{\alpha,\varepsilon}(x,y)f(y)\overline{g(x)}\,dw^+_{\alpha}(y) \, dw^+_{\alpha}(x). \label{12c}
\end{gather}
\end{proposition}
\begin{proof}
We follow the lines of the proof of \cite[Proposition 4.2]{ST2}, see also \cite[Proposition 3.2]{ST1}.
By Parseval's identity and \eqref{brs},
\begin{gather}
\langle \mathcal L^{-i\gamma}_{\alpha,\varepsilon,+}f,g\rangle_{L^2(\R_+,w^+_{\alpha})} =
\sum_{n\in\mathcal{N}_{\varepsilon}}
(2|n|+2|\alpha|+2d)^{-i\gamma} \, \langle f, h_n^{\alpha}\rangle_{L^2(\R_+,w^+_{\alpha})}\,
	 \langle h_{n}^{\alpha}, g \rangle_{L^2(\R_+,w^+_{\alpha})}.
\label{13c}
\end{gather}
To f\/inish the proof it is now suf\/f\/icient to verify that the right-hand sides of \eqref{12c}
and \eqref{13c} coincide. This task means justifying the possibility of changing the order
of integration, summation and dif\/ferentiation in the relevant expressions, see the proof of
Proposition 4.2 in \cite{ST2}. The details are rather elementary and thus are omitted.
The key estimate
\begin{gather*}
\int_{\R_+}\int_{\R_+}\int_0^\infty \big|\partial_t G^{\alpha,\varepsilon}_t(x,y)
\big|\,dt\, |\overline{g(x)}f(y)|\,dy\,dx<\infty
\end{gather*}
is easily verif\/ied by means of Lemma \ref{der_est}.
Another important ingredient (implicit in the proof of \cite[Proposition 4.2]{ST2}) is a suitable estimate
for the growth of the underlying eigenfunctions. In the present setting it is suf\/f\/icient to know that
\begin{gather*}
|h_n^{\alpha}(x)| \lesssim \prod_{i=1}^d \Phi_{n_i}^{\alpha_i}(x_i), \qquad x \in \mathbb{R}^d_+,
\end{gather*}
where
\begin{gather*}
\Phi_{n_i}^{\alpha _i} (x_i) = x_i^{-\alpha_i-1\slash 2} \left\{ \begin{array}{ll}
1,  &  0<x_i \le 4(n_i+\alpha _i+1); \\
\exp(- c x_i),  &  x_i > 4(n_i+\alpha _i+1). \end{array} \right.
\end{gather*}
This follows from Muckenhoupt's generalization \cite{Mu2} of the classical
estimates due to Askey and Wainger \cite{AW}.
\end{proof}

The theorem below says that the kernel
$K_\gamma^{\alpha,\varepsilon}(x,y)$ satisf\/ies standard estimates in the sense of the
homogeneous space $(\mathbb{R}^d_+,w^+_{\alpha},\|\cdot\|)$.
The corresponding proof is located in Section \ref{sec:kernel} below.
Denote $B^+(x,r) = B(x,r)\cap \mathbb{R}^d_+$.
\begin{theorem} \label{ker_est}
Given $\alpha \in [-1\slash 2, \infty)^d$ and $\varepsilon \in \mathbb{Z}_2^d$,
the kernel $K_\gamma^{\alpha,\varepsilon}(x,y)$ satisfies the growth condition
\begin{gather*}
|K_\gamma^{\alpha,\varepsilon}(x,y)| \lesssim
\frac{1}{w^+_{\alpha}(B^+(x,\|y-x\|))}, \qquad x,y \in \R_+, \quad x \neq y,
\end{gather*}
and the smoothness condition
\begin{gather*}
\|\nabla_{\! x,y} K_\gamma^{\alpha,\varepsilon}(x,y)\| \lesssim \frac{1}{\|x-y\|}
	\frac{1}{w^+_{\alpha}(B^+(x,\|y-x\|))}, \qquad x,y \in \R_+, \quad x \neq y.
\end{gather*}
\end{theorem}

From Theorem \ref{ker_est} and Proposition \ref{sto} we conclude that
$\mathcal L^{-i\gamma}_{\alpha,\varepsilon,+}$ is a Calder\'on--Zygmund operator.
Thus Lemma \ref{lem:main} follows from the general theory, see \cite{Duo}.

\begin{remark}
The results of this section can be generalized in a straightforward manner by
considering weighted $L^p$ spaces. By the general theory, each
$\mathcal L^{-i\gamma}_{\alpha,\varepsilon,+}$ extends to a bounded operator
on $L^p(\R_+,Wdw^+_{\alpha})$, $W\in A_p^{\alpha}$, $1<p<\infty$,
and to a bounded operator from $L^1(\R_+,Wdw^+_{\alpha})$ to $L^{1,\infty}(\R_+,Wdw^+_{\alpha})$,
$W \in A_1^{\alpha}$; here $A_p^{\alpha}$ stands for the Muckenhoupt class of $A_p$ weights associated
with the space $(\R_+,w^+_{\alpha},\|\cdot\|)$.
Consequently, analogous mapping properties hold for $\mathcal L^{-i\gamma}_{\alpha}$,
with ref\/lection invariant weights satisfying $A^{\alpha}_p$ conditions when restricted to $\R_+$
(or, equivalently, satisfying $A_p$ conditions related to the whole space $(\R,w_{\alpha},\|\cdot\|)$).
\end{remark}

\begin{remark}
With the particular $\varepsilon_0=(0,\ldots,0)$ the operator
$\mathcal L^{-i\gamma}_{\alpha,\varepsilon_0,+}$ coincides, up to a constant factor,
with the same imaginary power of the Laguerre Laplacian investigated
in \cite{NS1}. Therefore the results of this section deliver also analogous results
in the setting of \cite{NS1}.
\end{remark}

\section{Kernel estimates } \label{sec:kernel}

This section is mainly devoted to the proof of the standard estimates stated in Theorem \ref{ker_est}.
The proof follows the pattern of the proof of Proposition 3.1 in \cite{NS1}, see also \cite{NS2}.
We use the formula
\begin{gather*}
G^{\alpha,\varepsilon}_t(x,y) =\frac{1}{2^d}
	\Big( \frac{1-\zeta^2}{2\zeta}\Big)^{d + |\alpha| + |\varepsilon|}
		(xy)^{\varepsilon} \!\int_{[-1,1]^d}\!\!
    \exp\Big(\!-\frac{1}{4\zeta} q_{+}(x,y,s) - \frac{\zeta}{4} q_{-}(x,y,s)\!\Big)
    	  \Pi_{\alpha+\varepsilon}(ds),
\end{gather*}
where
\begin{gather*}
q_{\pm}(x,y,s) = \|x\|^2 + \|y\|^2 \pm 2 \sum_{i=1}^d x_i y_i s_i
\end{gather*}
(for the sake of brevity we shall often write shortly $q_{+}$ or $q_{-}$ omitting the arguments) and
$t\in(0,\infty)$ and $\zeta \in (0,1)$ are related by $\zeta = \tanh t$, so that
\begin{gather} \label{tzeta}
t = t(\zeta) = \frac{1}{2} \log \frac{1+ \zeta}{1-\zeta};
\end{gather}
eventually, $\Pi_{\alpha}$ denotes the product measure $\bigotimes\limits_{i=1}^d \Pi_{\alpha_i}$,
where $\Pi_{\alpha_i}$ is determined by the density
\begin{gather*}
\Pi_{\alpha_i}(ds) = \frac{(1-s^2)^{\alpha_i-1\slash 2}ds}{\sqrt{\pi}
	2^{\alpha_i}\Gamma{(\alpha_i+1\slash 2)}}, \qquad s \in (-1,1),
\end{gather*}
when $\alpha_i > -1\slash 2$, and in the limiting case of $\alpha_i = -1\slash 2$,{\samepage
\begin{gather*}
\Pi_{-1\slash 2} = \frac{1}{\sqrt{2\pi}} \big( \eta_{-1} + \eta_1 \big)
\end{gather*}
($\eta_{-1}$ and $\eta_1$ denote point masses at $-1$ and $1$, respectively).}

By the change of variable \eqref{tzeta} the kernels \eqref{kernel} can be expressed as
\begin{gather} \label{k_r}
K_\gamma^{\alpha,\varepsilon}(x,y) = \int_{[-1,1]^d} \Pi_{\alpha + \varepsilon}(ds) \int_0^1
	\beta_{d,\alpha+\varepsilon}(\zeta) \psi_{\zeta}^{\varepsilon}(x,y,s) \, d\zeta,
\end{gather}
where
\begin{gather*}
\psi_{\zeta}^{\varepsilon}(x,y,s) = (xy)^{\varepsilon}
	\exp\left(-\frac{1}{4\zeta} q_{+}(x,y,s) - \frac{\zeta}{4} q_{-}(x,y,s)\right)
\end{gather*}
and
\begin{gather*}
\beta_{d,\alpha}(\zeta) = \frac{{2^{1-d-i\gamma}}}{\Gamma(i\gamma)} \left( \frac{1-\zeta^2}{2\zeta}
	 \right)^{d+|\alpha|} \frac{1}{1-\zeta^2} \left( \log \frac{1+\zeta}{1-\zeta} \right)^{i\gamma -1}.
\end{gather*}
Notice that $|\beta_{d,\alpha}(\zeta)|$ coincides, up to a constant factor, with
$\beta^0_{d,\alpha}(\zeta)$ def\/ined in \cite[(5.4)]{NS1}.

The application of Fubini's theorem that was necessary to get \eqref{k_r} is also justif\/ied since,
in fact, the proof (to be given below) of the f\/irst estimate in Theorem \ref{ker_est} contains
the proof of
\begin{gather*}
\int_{[-1,1]^d} \Pi_{\alpha + \varepsilon}(ds) \int_0^1
	\big|\beta_{d,\alpha+\varepsilon}(\zeta) \psi_{\zeta}^{\varepsilon}(x,y,s)\big| \, d\zeta < \infty,
	\qquad x \neq y.
\end{gather*}

For proving Theorem \ref{ker_est} we  need a specif\/ied version of
\cite[Corollary 5.2]{NS1} and a slight extension of \cite[Lemma 5.5 (b)]{NS1}
(the proof of the latter result in \cite{NS1} is given under assumption $k \ge 1$, but in fact
it is also valid for any real $k$ provided that the constant factor in the def\/inition of
$\beta^k_{d,\alpha}(\zeta)$ is neglected; in particular, $k=0$ can be admitted).
\begin{lemma}
\label{m_lem}
Assume that $\alpha \in [-1\slash 2,\infty)^d$. Let $b \ge 0$ and $c>0$ be fixed. Then, we have
\begin{gather*}
  \emph{(a)} \quad
	\big( |x_1 \pm y_1 s_1| + |y_1 \pm x_1 s_1| \big) \exp\left( - c\frac{1}{\zeta} q_{\pm}(x,y,s) \right)
     \lesssim \zeta^{\pm 1\slash 2},\\
  \emph{(b)} \quad
\int_0^1 \big| \beta_{d,\alpha}(\zeta)\big| \zeta^{-b}
    \exp\left( -c \frac{1}{\zeta} q_{+}(x,y,s) \right)   d\zeta
        \lesssim \big(q_{+}(x,y,s)\big)^{-d-|\alpha|-b},
\end{gather*}
uniformly in $x,y \in \R_+$, $s \in [-1,1]^d$, and also in $\zeta \in (0,1)$ if (a) is considered.
\end{lemma}

We also need the following generalization of \cite[Proposition 5.9] {NS1}, cf.~\cite[Lemma 5.3]{NS2}.

\begin{lemma} \label{lemhom}
Assume that $\alpha \in [-1\slash 2, \infty)^d$ and let $\delta,\kappa \in [0,\infty)^d$ be fixed.
Then for $x,y\in\R_+$, $x \neq y$,
\begin{gather*}
(x+y)^{2\delta} \int_{[-1,1]^d}\Pi_{\alpha+\delta+\kappa}(ds)
	\; \big(q_{+}(x,y,s)\big)^{-d - |\alpha| -|\delta|} \lesssim \frac{1}{w^+_{\alpha}(B^+(x,\|y-x\|))}
\end{gather*}
and
\begin{gather*}
 (x+y)^{2\delta}\int_{[-1,1]^d} \Pi_{\alpha+\delta+\kappa}(ds)\;
    \big(q_{+}(x,y,s)\big)^{-d - |\alpha| -|\delta| - 1\slash 2} \lesssim \frac{1}{\|x-y\|} \;
 \frac{1}{w^+_{\alpha}(B^+(x,\|y-x\|))}.
\end{gather*}
\end{lemma}

\begin{proof}[Proof of Theorem \ref{ker_est}]
The growth estimate is rather straightforward.
Using Lemma \ref{m_lem} (b) with $b=0$ and observing that $(xy)^{\varepsilon} \le (x+y)^{2\varepsilon}$
gives
\begin{gather*}
|K_\gamma^{\alpha,\varepsilon}(x,y)| \lesssim (x+y)^{2\varepsilon}
	\int_{[-1,1]^d} \Pi_{\alpha+\varepsilon}(ds) (q_+)^{-d-|\alpha|-|\varepsilon|}.
\end{gather*}
Now Lemma \ref{lemhom}, taken with $\delta = \varepsilon$ and $\kappa = (0,\ldots,0)$,
provides the desired bound.

It remains to prove the smoothness estimate. Notice that by symmetry reasons it is enough to show that
\begin{gather*}
|\partial_{x_1}K_\gamma^{\alpha,\varepsilon}(x,y)| + |\partial_{y_1}K_\gamma^{\alpha,\varepsilon}(x,y)|
\lesssim \frac{1}{\|x-y\|} \frac{1}{w^+_{\alpha}(B^+(x,\|y-x\|))}, \qquad x,y \in \R_+, \quad x \neq y.
\end{gather*}
Moreover, we can focus on estimating the $x_1$-derivative only. This is because in the f\/inal stroke
we shall use Lemma \ref{lemhom}, where the left-hand sides are symmetric in $x$ and $y$.
Thus we are reduced to estimating the quantity
\begin{gather*}
\mathcal{J} = \int_{[-1,1]^d} \Pi_{\alpha + \varepsilon}(ds) \int_0^1
	\big|\beta_{d,\alpha+\varepsilon}(\zeta) \partial_{x_1} \psi_{\zeta}^{\varepsilon}(x,y,s)\big| \, d\zeta
\end{gather*}
(passing with $\partial_{x_1}$ under the integral signs is legitimate, the justif\/ication being
implicitly contained in the estimates below, see the argument in \cite[pp.~671--672]{NS1}).

An elementary computation produces
\begin{gather*}
\partial_{x_1} \psi_{\zeta}^{\varepsilon}(x,y,s)  =
 \left[ (xy)^{\varepsilon} \left(- \frac{1}{2\zeta}(x_1+y_1 s_1)
	- \frac{\zeta}{2} (x_1-y_1 s_1) \right)+ \varepsilon_1
		y_1 (xy)^{\varepsilon-e_1} \right] \\
\phantom{\partial_{x_1} \psi_{\zeta}^{\varepsilon}(x,y,s)  =}{}
\times\exp\left( - \frac{1}{4\zeta} q_+ - \frac{\zeta}{4} q_{-} \right).
\end{gather*}
Hence, by Lemma \ref{m_lem} (a), we have
\begin{gather*}
 \big| \partial_{x_1} \psi_{\zeta}^{\varepsilon}(x,y,s) \big| \\
\qquad {} \lesssim
	(xy)^{\varepsilon} (\zeta^{-1\slash 2} + \zeta^{1\slash 2})
		\exp\left( - \frac{1}{8\zeta} q_+ - \frac{\zeta}{8} q_{-} \right)
	+ \varepsilon_1 y_1 (xy)^{\varepsilon -e_1}
	\exp\left( - \frac{1}{4\zeta} q_+ - \frac{\zeta}{4} q_{-} \right) \\
\qquad {} \lesssim
	(x+y)^{2\varepsilon} \zeta^{-1\slash 2} \exp\left( - \frac{1}{8\zeta} q_+  \right)
	+ \varepsilon_1 (x+y)^{2(\varepsilon - e_1\slash 2)}
	\exp\left( - \frac{1}{4\zeta} q_+  \right)
\end{gather*}
(notice that the second term above vanishes when $\varepsilon_1 = 0$).
Consequently,
\begin{gather*}
\mathcal{J} \lesssim   (x+y)^{2\varepsilon} \int_{[-1,1]^d} \Pi_{\alpha+\varepsilon}(ds)
	\int_0^1 \big| \beta_{d,\alpha+\varepsilon}(\zeta) \big| \zeta^{-1\slash 2}
	\exp\left( - \frac{1}{8\zeta} q_+ \right)  d\zeta \\
\phantom{\mathcal{J} \lesssim}{}  + \varepsilon_1 (x+y)^{2(\varepsilon - e_1\slash 2)} \int_{[-1,1]^d} \Pi_{\alpha+\varepsilon}(ds)
	\int_0^1 \big| \beta_{d,\alpha+\varepsilon}(\zeta) \big|
	\exp\left( - \frac{1}{4\zeta} q_+ \right) d\zeta.
\end{gather*}
Now, applying Lemma \ref{m_lem} (b) with either $b=1\slash 2$ or $b=0$ leads to
\begin{gather*}
\mathcal{J} \lesssim   (x+y)^{2\varepsilon} \int_{[-1,1]^d} \Pi_{\alpha+\varepsilon}(ds)
	(q_+)^{-d-|\alpha|-|\varepsilon|-1\slash 2} \\
\phantom{\mathcal{J} \lesssim}{}  + \varepsilon_1 (x+y)^{2(\varepsilon - e_1\slash 2)} \int_{[-1,1]^d} \Pi_{\alpha+\varepsilon}(ds)
	(q_+)^{-d-|\alpha|-|\varepsilon|}.
\end{gather*}
Finally, Lemma \ref{lemhom} with either $\delta=\varepsilon$ and $\kappa = (0,\ldots,0)$ or
(in case $\varepsilon_1 =1$) $\delta = \varepsilon - e_1\slash 2$ and $\kappa = e_1\slash 2$
delivers the required smoothness bound for $\mathcal{J}$.

The proof of Theorem \ref{ker_est} is complete.
\end{proof}

\begin{proof}[Proof of Lemma \ref{der_est}]
Recall that
\begin{gather*}
G^{\alpha,\varepsilon}_t(x,y) =\frac{1}{2^d}
	\left( \frac{1-\zeta^2}{2\zeta}\right)^{d + |\alpha| + |\varepsilon|}
		(xy)^{\varepsilon} \int_{[-1,1]^d} \Pi_{\alpha+\varepsilon}(ds) \,
    \exp\left(-\frac{1}{4\zeta} q_{+} - \frac{\zeta}{4} q_{-}\right),
\end{gather*}
where $t$ and $\zeta$ are related by $\zeta = \tanh t$.
Since $\zeta \in (0,1)$ and $\|x-y\|^2 \le q_{\pm} \le \|x+y\|^2$, we see that
\begin{gather*}
G^{\alpha,\varepsilon}_t(x,y) \lesssim
	\left( \frac{1-\zeta}{\zeta}\right)^{d + |\alpha| + |\varepsilon|} (xy)^{\varepsilon}
	\exp\left( -\frac{1}{4\zeta}\|x-y\|^2\right).
\end{gather*}
From this estimate the rapid decay easily follows.

To verify \eqref{g_e} we need f\/irst to compute $\partial_t G_t^{\alpha,\varepsilon}(x,y)$.
We get
\begin{gather*}
\partial_t G_t^{\alpha,\varepsilon}(x,y) =   \big(d+|\alpha|+|\varepsilon|\big) \frac{1+\zeta^2}{\zeta}
	G_t^{\alpha,\varepsilon}(x,y)
	 + \frac{1}{2^d} \left( \frac{1-\zeta^2}{2\zeta}\right)^{d + |\alpha| + |\varepsilon|} (1-\zeta^2)
	(xy)^{\varepsilon} \\
\phantom{\partial_t G_t^{\alpha,\varepsilon}(x,y) =}{}  \times \int_{[-1,1]^d} \Pi_{\alpha+\varepsilon}(ds)
	\left(\frac{1}{4\zeta^2} q_{+} - \frac{1}{4} q_{-}\right)
    \exp\left(-\frac{1}{4\zeta} q_{+} - \frac{\zeta}{4} q_{-}\right)
\end{gather*}
(here passing with $\partial_t$ under the integral can be easily justif\/ied).
Consequently, taking into account the estimates above,
\begin{gather*}
\big| \partial_t G_t^{\alpha,\varepsilon}(x,y) \big| \lesssim   \frac{1}{\zeta}
	\left( \frac{1-\zeta}{\zeta}\right)^{d + |\alpha| + |\varepsilon|} (xy)^{\varepsilon}
	\exp\left( -\frac{1}{4\zeta}\|x-y\|^2\right) \\
\phantom{\big| \partial_t G_t^{\alpha,\varepsilon}(x,y) \big| \lesssim }{}  + (1-\zeta)\left( \frac{1-\zeta}{\zeta}\right)^{d + |\alpha| + |\varepsilon|}
	(xy)^{\varepsilon} \frac{\|x+y\|^2}{\zeta^2} \exp\left( -\frac{1}{4\zeta}\|x-y\|^2\right).
\end{gather*}
This implies
\begin{gather*}
\int_0^{\infty} \big| \partial_t G_t^{\alpha,\varepsilon}(x,y)\big| \, dt =
	\int_0^1 \left| \partial_t G_t^{\alpha,\varepsilon}(x,y)\big|_{t = \tanh^{-1} \zeta} \right|
		\frac{d\zeta}{1-\zeta^2} \\
\phantom{\int_0^{\infty} \big| \partial_t G_t^{\alpha,\varepsilon}(x,y)\big| \, dt}{}
\lesssim      (xy)^{\varepsilon} \int_0^1
			({1-\zeta})^{d + |\alpha| + |\varepsilon|-1} \zeta^{-(d + |\alpha| + |\varepsilon|-1)}
			\exp\left( -\frac{1}{4\zeta}\|x-y\|^2\right)   d\zeta \\
\phantom{\int_0^{\infty} \big| \partial_t G_t^{\alpha,\varepsilon}(x,y)\big| \, dt=}{}
+ (xy)^{\varepsilon} \|x+y\|^2 \int_0^1  \zeta^{-(d + |\alpha| + |\varepsilon|-2)}
			\exp\left( -\frac{1}{4\zeta}\|x-y\|^2\right)   d\zeta.
\end{gather*}
Now using the fact that $d + |\alpha| + |\varepsilon|>0$ and
$\sup_{u>0} u^a \exp(-Au)< \infty$ for any f\/ixed $A>0$ and $a \ge 0$, leads to the bound
\begin{gather*}
\int_0^{\infty} \big| \partial_t G_t^{\alpha,\varepsilon}(x,y)\big| \, dt \lesssim
	\frac{(xy)^{\varepsilon}}{\|x-y\|^{2(d + |\alpha| + |\varepsilon|+1)}} +
	\frac{(xy)^{\varepsilon} \|x+y\|^2}{\|x-y\|^{2(d + |\alpha| + |\varepsilon|+2)}}.
\end{gather*}
The conclusion follows.
\end{proof}

\pdfbookmark[1]{References}{ref}
 \LastPageEnding

\begin{thebibliography}{99}

\footnotesize\itemsep=0pt

\bibitem{AW}
 Askey R., Wainger S.,
 Mean convergence of expansions in Laguerre and Hermite series,
\emph{Amer. J. Math.} \textbf{87} (1965), 695--708.

\bibitem{Du2}
 Dunkl C.F.,
Ref\/lection groups and orthogonal polynomials on the sphere,
\emph{Math. Z.} \textbf{197} (1988), 33--60.

\bibitem{Du1}
 Dunkl C.F.,
{Dif\/ferential-dif\/ference operators associated to ref\/lection groups},
\emph{Trans. Amer. Math. Soc.} \textbf{311} (1989), 167--183.

\bibitem{Duo}
 Duoandikoetxea J.,
Fourier analysis,
{\it Graduate Studies in Mathematics}, Vol.~29, American Mathematical Society, Providence, RI, 2001.

\bibitem{Leb}
 Lebedev N.N.,
Special functions and their applications,
Dover Publications, Inc., New York, 1972.

\bibitem{M}
 Muckenhoupt B.,
On certain singular integrals,
\emph{Pacific J. Math.} \textbf{10} (1960), 239--261.

\bibitem{Mu2}
Muckenhoupt B.,
{Mean convergence of Hermite and Laguerre series.~II},
\emph{Trans. Amer. Math. Soc.} \textbf{147} (1970), 433--460.

\bibitem{NS1}
 Nowak A., Stempak K.,
Riesz transforms for multi-dimensional Laguerre function expansions,
\emph{Adv. Math.} \textbf{215} (2007), 642--678.

\bibitem{NS2}
 Nowak A., Stempak K.,
Riesz transforms  for  the Dunkl harmonic oscillator,
\emph{Math. Z.}, to appear, \href{http://arxiv.org/abs/0802.0474}{arXiv:0802.0474}.

\bibitem{Ros}
 Rosenblum M.,
Generalized Hermite polynomials and the Bose-like oscillator calculus,
in Nonselfadjoint Operators and Related Topics (Beer Sheva, 1992),
\emph{Oper. Theory Adv. Appl.}, Vol.~73, Birkh\"auser, Basel, 1994, 369--396.

\bibitem{R1}
 R\"osler M.,
{Generalized Hermite polynomials and the heat equation for Dunkl operators},
\emph{Comm. Math. Phys.} \textbf{192} (1998), 519--542, \href{http://arxiv.org/abs/q-alg/9703006}{q-alg/9703006}.

\bibitem{R2}
 R\"osler M.,
{One-parameter semigroups related to abstract quantum models of Calogero types},
in Inf\/inite Dimensional Harmonic Analysis (Kioto, 1999),
 Gr\"abner, Altendorf, 2000, 290--305.

\bibitem{R3}
 R\"osler M.,
{Dunkl operators: theory and applications},
in Orthogonal Polynomials and Special Functions (Leuven, 2002),
 \emph{Lecture Notes in Math.}, Vol.~1817, Springer, Berlin, 2003, 93--135, \href{http://arxiv.org/abs/math.CA/0210366}{math.CA/0210366}.

\bibitem{ST1}
 Stempak K., Torrea J.L.,
{Poisson integrals and Riesz transforms for Hermite function expansions with weights},
\emph{J. Funct. Anal.} \textbf{202} (2003), 443--472.

\bibitem{ST2}
 Stempak K., Torrea J.L.,
Higher Riesz transforms and imaginary powers associated to the harmonic oscillator,
\emph{Acta Math. Hungar.} \textbf{111} (2006), 43--64.

\end{thebibliography}
\end{document}